\documentclass[12pt]{amsart}
\usepackage{amsfonts}
\usepackage{amsmath}
\usepackage{graphicx}
\usepackage[ansinew]{inputenc}

\newtheorem{teo}{Theorem}[section]
\newtheorem{prop}{Proposition}[section]

\newtheorem{lema}{Lemma}[section]
\newtheorem{defi}{Definition}[section]

\newtheorem{nota}{Remark}[section]

\def \e{\hskip 10pt }
\def \ecu{\hskip 5pt }

\def \Dem{{\bf Proof. }}

\def \ep{\epsilon}
\def \sign{\mbox{sign }\! \!}

\def \bbbr{\mathbb R}

\def \bbbn{\mathbb N}
\newcommand{\re}{I \! \! R}

\def \num{${\mbox n}^{\underline{\mbox{\tiny o}}}$ }
\def \12{\frac{1}{2}}

\def \hate {\hskip 2pt \hat{} \hskip 2pt}

\def \al {\alpha }

\def \S{{\mathcal S}}
\def \T{{\mathcal T}}

\def \MB{{\mathcal MB}}

\def \t2MB{\tilde{\tilde{{\mathcal MB}}}}

\def \cqd{\vskip 20pt}

\def\lpq{{L^{p,q}}}

\def\ds{\displaystyle}

\title{Bilinear multipliers on Lorentz spaces}%

\author{Francisco Villarroya}
\address{Departament d'An\`alisi Matem\`atica, Universitat de Val\`encia,
46100 Burjassot, Val\`encia, Spain}
\email{paco.villarroya@uv.es}
\thanks{The author has been partially supported by grants DGESIC
PB98-1246 and BMF 2002-04013}


\begin{document}
\maketitle

\begin{abstract}
%Let $m$ a measurable function defined on $\bbbr^{2}$. A
%bilinear multiplier defined by $m$ is the operator
%that acts over Schwartz functions in the following way:
%$
%B_{m}(f,g)(x)=\int_{\re^{2}}\hat{f}(\xi )\hat{g}(\eta )m(\xi ,\eta )
%e^{2\pi i(\xi +\eta )x}d\xi d\eta
%$.
We give one sufficient and two necessary conditions for boundedness
between Lebesgue or Lorentz spaces of several classes of bilinear
multiplier operators closely
connected with the bilinear Hilbert transform.
\end{abstract}

%%%%%%%%%%%%%%%%%%%%%%%%%%%%%%%%%%%%%%%%%%%%%%%%%%%%%%%%%
%%%%%%%%%%%%%%%%%%%%%%%%%%%%%%%%%%%%%%%%%%%%%%%%%%%%%%%%%

\section{Introduction.}
The bilinear Hilbert transform with parameter
$\alpha \in \bbbr $ is the operator given by
$$
H_{\alpha }(f,g)(x)=\frac{1}{\pi} p.v.\int f(x-t)g(x-\alpha t)\frac{dt}{t}
$$
initially defined for functions in the Schwartz class.
Notice that $H_{0}(f,g)=H(f)g $ and $H_{1}(f,g)=H(fg)$ where $H(f)$ is the
classical Hilbert transform. So
$H_{\alpha }$ can be seen as an intermediate step between both operators.

The bilinear Hilbert transform has been extensively studied since 1965 when A.
Calder\'on set the hypothesis of its boundedness from $L^{2}\times L^{\infty }$
into $L^{2}$ while he was working on the Hilbert transform defined over
Lipschitz curves (see \cite{cal}). After several years of research and using original ideas of
C. Fefferman \cite{feff}, M. Lacey and C. Thiele finally answered this
question when they proved the following

\begin{teo}\label{laceythiele}
 For each triple $(p_1,p_2,p_3)$ such
that $1<p_1,p_2\le\infty$,
$1/p_1+1/p_2=1/p_3$ and
$p_3>2/3$ and each $\alpha\in \bbbr \setminus\{ 0,1\}$ there exists
$C(\alpha,p_1,p_2)>0$ for which
$$\|H_{\alpha}(f,g)\|_{p_3}\le
C(\alpha,p_1,p_2)\|f\|_{p_1}\|g\|_{p_2}
$$
for all $f,g$ in the Schwartz class
\end{teo}
\noindent in two papers (\cite{lath1},
\cite{lath2}) published in 1997 and 1999 respectively.
See also \cite{thielehab} for a unified proof.

From then a great deal of generalizations and extensions of
this seminal work have appeared such that: \cite{GN1}, \cite{GN2} and
\cite{MTT} related to the modification
of the kernel of the operator, \cite{grali1} related to
uniform estimates in the same inequality, \cite{lacey2} related to maximal
results, \cite{mutath2} to uniform estimates with generalized kernels.
%and \cite{grataoter,grataoter} to operators that act over functions of several
%variables. All these papers extend the techniques used in
%\cite{feff}, the time-frequency analysis, which basically consists on
%decompositions of functions well located both in time and frequency.

The present paper shows two sufficient and one necessary conditions
for boundedness of different types of bilinear multipliers some of which include the
bilinear Hilbert transform.

\section{Preliminaries, notation and definitions.}
 Given a measurable function $f$ we denote its
distribution function by $m_f(\lambda) = m(\{ x\in \bbbr : |f(x)| > \lambda\} )$ and
its nonincreasing rearrangement by $f^*(t) = \inf\{\lambda>0 : m_f(\lambda) \leq t\}$.
The Lorentz space $\lpq$ consists of those
measurable functions $f$ such that $\|f\|^*_{p,q} < \infty$, where
$$
\|f\|^*_{p,q} = \left\{
\begin{array}{ll}
\ds{\left\{ \frac{q}{p} \int_0^\infty t^{\frac{q}{p}} f^*(t)^q \frac{dt}{t}\right\}^\frac{1}{q},} &
0<p<\infty ,\ 0 < q < \infty , \\
\ds{\sup_{t> 0} t^{\frac{1}{p}} f^*(t) }  & 0 < p \leq \infty , \ q=\infty .
\end{array} \right.
$$

Reader is referred to \cite{BS}
%and \cite{H}
for basic information on Lorentz spaces.

The interpolation result we are going to use is a trilinear version of
Riesz-Thorin interpolation theorem over tuplas of spaces. Since we will
use it with positive integral operators
$$
\int_{\re }f(x-t)g(x-\alpha t)K(t)dt
$$
where $K$ is a positive function, we state
the theorem in this setting.

\begin{teo}\label{interpolate}
Let $0<p_{i,j}\leq \infty $ for $i=1,\ldots ,n$, $j=0,1,2,3$. Let $T$ a
positive trilinear integral operator such that
$
T:L^{p_{i,0}}\times L^{p_{i,1}}\times L^{p_{i,2}} \rightarrow L^{p_{i,3}}
$
is bounded for $i=1,\ldots ,n$
with $\| T\|_{i}\leq M_{i}$.

Then
$
T:L^{p_{0}}\times L^{p_{1}}\times L^{p_{2}} \rightarrow L^{p_{3}}
$
is bounded
for $\frac{1}{p_{j}}=\sum_{i=1}^{n}\frac{\theta_{i}}{p_{i,j}}$, for
$j=0,1,2,3$ where $0\leq \theta_{i} \leq 1$ and
$\sum_{i=1}^{n}\theta_{i}=1$. Moreover,
$\| T\|\leq \prod_{i=1}^{n}M_{i}^{\theta_{i} }$.
\end{teo}

A proof of this theorem between pair of spaces can be seen in \cite{BS}
page 185 for the linear case and 202 for the multilinear case. The
extension to tuplas of spaces is trivial from that result.

%\vskip 10pt
%We make now explicit the kind of norm we use in the Schwartz class in order to
%work with temperated distributions.
%\begin{defi}{\bf (Topology $\T_{\S }$).}
%The Schwartz class, $\S (\bbbr^{n})$, turns out to be a Fr\'echet space with
%the topology
%$\T_{\S }$ defined by the family of seminorms
%$$
%\rho_{m}(f)=\sup_{|\alpha |,|\beta |\leq m}
%\| f_{\al ,\be }\|_{\infty }
%$$
%where $\alpha, \beta \in \bbbn^{n}$, $m\in \bbbn $  and $f_{\al ,\be }(x)=x^{\al }D^{\be }f(x)$.

%Another way of defining the same topology is by saying that
%a net $\{ \varphi_{j}\}_{j\in J}\subset \S (\bbbr^{n})$ converges to
%$\varphi $ if for all $\alpha ,\beta \in \bbbn^{n}$
%we have that $\{ x^{\alpha }D^{\beta}\varphi_{j}(x)\}_{j\in J}$
%converges uniformly to $x^{\alpha }D^{\beta}\varphi (x)$ .

%Obviously, we call temperate distribution to every element of the dual space
%of $(\S (\bbbr^{n}), \T_{\S })$.
%\end{defi}

We set some frequently used notation.
For every $x,y\in \bbbr $ we denote the translation operator by $T_{y}f(x)=f(x-y)$ and
the modulation operator by $M_{y}f(x)=f(x)e^{2\pi iyx}$
while for all $p\in \bbbr $ and $t\neq 0$ we denote the dilation operators by
$D_{t}^{p}f(x)=t^{-\frac{1}{p}}f(t^{-1}x)$ and
$D_{t}f(x)=D_{t}^{\infty }f(x)=f(t^{-1}x)$.
These operators show certain symmetries when the Fourier transform acts over
them. In particular, the transform of a translation is a modulation,
$(T_{y}f)\hate{}=M_{-y}\hat{f}$, the transform of a modulation is a translation,
$(M_{y}f)\hate{}=T_{y}\hat{f}$
and the transform of a dilation is its dual dilation,
$(D_{t}^{p}f)\hate{}=\sign(t)\ecu D_{t^{-1}}^{p'}\hat{f}$.

For the dilation operator we trivially have that
$
\| D_{t}^{r}f\|_{p,q}=|t|^{\frac{1}{p}-\frac{1}{r}}\| f\|_{p,q}
$.
Sometimes we will also use the notation $K_{\ep }$ for the change of scale
normalized to the $L^{1}$ norm, that is,
$
K_{\ep }(x)=\ep^{-1}K(\ep^{-1}x)=D_{\ep }^{1}K(x)
$.

%We see now the definition of the bilinear operator we are going to work with.
%The most general way of define operators like the bilinear Hilbert transform
%is by working distributionally. In some sense, these operators can be seen as
%modifications or generalizations of convolution operators and thus, they can
%be defined in the same way the convolution of a distribution and a function
%can be done, that is, functionally and distributionally. We will work only
%with the functional one.

The bilinear operators we are going to work with can be seen as generalizations of convolution operators. Thus,
as in the case of the convolution of a distribution and a function,
they can be defined functionally and distributionally. We will work only
with the functional definition.

\begin{defi}
Let $u$ be a distribution. For every
$\alpha \in \bbbr $ and every
$f,g\in C_{0}^{\infty }$ we define the function
$$
H_{u,\alpha }(f,g)(x)=(u,D_{-1}T_{-x}f\cdot D_{-\al^{-1}}T_{-x}g)
$$
for all $x\in \bbbr $. We will say that $H_{u,\alpha }$ is
a generalized bilinear Hilbert transform associated to $u$ and $\alpha $ or
just a BHT for short.
\end{defi}

In this way, if $K$ is a locally integrable function for instance this
definition leads to the expression
\begin{equation}\label{nucleok}
H_{K,\alpha }(f,g)(x)=\int_{\re}f(x-t)g(x-\alpha t)K(t)dt
\end{equation}
which is well defined for all $\alpha, x\in \bbbr $ and for every
$f,g$ bounded functions such that at least one of them
has compact support if $\al \neq 0$ or $f$ has
compact support if $\al =0$.

We give the following
\begin{defi}Let $\al\in \bbbr $ and $u$ be a distribution. Let
$0<p_{i}<\infty $, $0<q_{i}\leq \infty $, $i=1,2,3$. We say that
$H_{u,\alpha }$ is
$(p_{i},q_{i})_{i=1,2,3}$ bounded if it can be extended to a
bounded operator from $L^{p_{1},q_{1}}\times L^{p_{2},q_{2}}$ into
$L^{p_{3},q_{3}}$. This is possible if there exists a constant
$C>0$ depending of $u$, $\alpha $ and $p_{i},q_{i}$ such that
$
\| H_{u,\alpha }(f,g)\|_{p_{3},q_{3}}\leq C\| f\|_{p_{1},q_{1}}
\| g\|_{p_{2},q_{2}}
$,
for all $f$ and $g$ in some appropriate dense subspaces.
%In particular if $p_{i}=q_{i}$ we say that $H_{u,\alpha }$ is
%$(p_{1},p_{2},p_{3})$ strong
%($
%\| H_{u,\alpha }(f,g)\|_{p_{3}}\leq C\| f\|_{p_{1}}\| g\|_{p_{2}}
%$)
%and if $p_{i}=q_{i}$ for $i=1,2$ and $q_{3}=\infty  $ we say that
%$H_{u,\alpha }$ is
%$(p_{1},p_{2},p_{3})$ weak
%($
%\| H_{u,\alpha }(f,g)\|_{p_{3},\infty }\leq C\| f\|_{p_{1}}\| g\|_{p_{2}}
%$).
\end{defi}

%Although it will not be used in the paper, one can
%define the operator distributionally in the following way
%
%\begin{defi}{\bf (distributional definition)}
%Let $u$ a distribution. Then, for all
%$\alpha \in \bbbr $ and for every
%$f,g\in C_{0}^{\infty }$ we can define the distribution
%$$
%(H_{u,\alpha }(f,g),h)=(u,H_{h,\alpha^{-1}}(D_{-1}f,D_{-\alpha^{-1}}g))
%$$
%for all $h\in C_{0}^{\infty }$ where
%$H_{h,\alpha^{-1}}(D_{-1}f,D_{-\alpha^{-1}}g)$ is defined like in
%(\ref{nucleok}). Notice that $H_{u,\alpha
%}(f,g)$ is in this way well defined once one see that if $f,g$
%have compact support, so do
%$H_{h,\alpha^{-1}}(D_{-1}f,D_{-\alpha^{-1}}g)$ (we only need $h$ to be
%locally integrable) and that if
%$f,g,h$ belong to $C^{\infty }$, so do
%$H_{h,\alpha^{-1}}(D_{-1}f,D_{-\alpha^{-1}}g)$.
%\end{defi}
%
%Moreover, it's very easy to see that both operators are equal
%distributionally.

In the same way that convolution and linear multiplier operators are intimately
related, so do are the operators previously defined and the following ones:
\begin{defi}
Let $m$ be a bounded measurable function in $\bbbr^{2}$. For every $x\in \bbbr $ and $f,g\in \S $ we define
the operator
$$
B_{m}(f,g)(x)
=\int_{\re^2}\hat{f}(\xi )\hat{g}(\eta )m(\xi ,\eta )e^{2\pi i(\xi +\eta )x}
d\xi d\eta
$$

Let $p_{i}>0$.
We say that $m$ is a $(p_{1},p_{2},p_{3})$ multiplier or just a bilinear
multiplier if the operator can be extended to a bounded operator from
$L^{p_{1}}\times L^{p_{2}}$ to $L^{p_{3}}$.
We denote by
$\| \cdot \|_{\MB_{p_{1},p_{2},p_{3}}}$ the minimum
constant that satisfy the inequality
$
\| B_{m}(f,g)\|_{p_{3}}\leq C \| f\|_{p_{1}}\| g\|_{p_{2}}
$
for all functions $f,g\in \S$.
%In this way,
%$\MB_{p_{1},p_{2},p_{3}}$ turns out to be a complete space.
\end{defi}

The relationship between both kind of operators is the following: if
$K$ is, we say, an integrable function then
$$
\int_{\re}f(x-t)g(x-\alpha t)K(t)dt
=\int_{\re^2 }\hat{f}(\xi )\hat{g}(\eta )
\widehat{K}(\xi +\alpha \eta)e^{2\pi i(\xi +\eta )x} d\xi d\eta
$$
and so, both operators can be regarded as generalization of convolution
operators or as generalization of linear multiplier operators.

We finally state several of their properties related to
invariance by traslation, commutativity and duality:
\begin{eqnarray}
\label{invtrasl}
H_{T_{y}u,\alpha }(f,g)&=&H_{u,\alpha }(T_{y}f,T_{\alpha y}g)\\
\label{conmu}
H_{u,\alpha }(f,g)&=&\sign(\alpha )H_{D_{\alpha }^{1}u,\alpha^{-1}}(g,f)\\
\label{dua}
\Big\langle h,H_{u,\alpha }(f,g)\Big\rangle &=& \Big\langle H_{D_{-1}u,1-\alpha
}(h,g),f\Big\rangle
\end{eqnarray}

\section{Three conditions for boundedness}
We introduce three results of boundedness which can be summarized as follows.
We first give a necessary condition obtained when we study the operator
acting over gaussian functions. Then we also give a sufficient
condition which is the
generalization of Young inequality to this class of non-convolution
operators. The
third one is another sufficient condition for the second class of operators we have defined.

\subsection{Gaussians looking for necessary conditions.}
We use the fact that the BHT over gaussian functions
has a particularly easy expresion in order to get
necessary conditions for its boundedness when the kernel is a temperate
distribution. We get in this way two conditions of boundedness: one over the
spaces between which the BHT can be bounded and another one over the kernel
itself. We work with Lorentz spaces just for the sake of generality.
We begin with a technical lemma.
\begin{lema}\label{aproxid}
Let $G\in \S $ such that $\hat{G}(0)=1$.  Let $(G_{\ep })_{\ep >0}$ an
approximate identity with
$G_{\ep }=D_{\ep }^{1}G$.
Then for all $\varphi \in \S $,
$(G_{\ep }*\varphi )_{\ep >0}$ converges to $\varphi $ in the topology
of the Schwartz class $\T_{\S }$.
\end{lema}
\Dem
We need to prove that for every $n,m\in \bbbn $,
$
\lim_{\ep \rightarrow 0^{+}}
\| (G_{\ep }*\varphi )_{n,m} - \varphi_{n,m}\|_{\infty }=0
$
where we define $\varphi_{n,m}(x)=x^{n}\varphi^{m)}(x)$.
If $c_{n,k}$
%$
%=\left( \begin{array}{c}
%n \\
%k
%\end{array}\right)
%$
denote the combinatorial number $n$ over $k$ then for
$x\in \bbbr $ and $\ep >0$ we have
$$
x^{n}(G_{\ep }*\varphi )^{m)}(x)=x^{n}(G_{\ep }*\varphi^{m)})(x)
=\int_{\re }(x-t+t)^{n}G_{\ep }(t)\varphi^{m)}(x-t)dt
$$
$$
=\sum_{k=0}^{n}c_{n,k}
\int_{\re }t^{k}D_{\ep }^{1}G(t)(x-t)^{n-k}\varphi^{m)}(x-t)dt
=\sum_{k=0}^{n}c_{n,k}\ep^{k}(D_{\ep }^{1}(G_{k,0})*\varphi_{n-k,m})(x)
$$
Thus,
$$
|(G_{\ep }*\varphi )_{n,m}(x) - \varphi_{n,m}(x)|
\leq |(G_{\ep }*\varphi_{n,m})(x)-\varphi_{n,m}(x)|
+\sum_{k=1}^{n}c_{n,k}\ep^{k}\| G_{k,0}\|_{1}\| \varphi_{n-k,m}\|_{\infty }
$$
and for $a=\max{(n,m)}$, $\rho_{r}(\varphi )=\sup_{m,n\leq r}\| \varphi_{n,m}\|_{\infty }$
$$
\| (G_{\ep }*\varphi )_{n,m}-\varphi_{n,m}\|_{\infty }
\leq \| G_{\ep }*\varphi_{n,m}-\varphi_{n,m}\|_{\infty }
+((\ep +1)^{n}-1)\max_{0\leq k\leq a}\| G_{k,0}\|_{1}\ecu \rho_{a}(\varphi )
$$
This proves the result by the main property of an approximate identity.

\begin{prop}\label{gaussian}
Let $\al <0$ and $p_{i},q_{i}>0$ for $i=1,2,3$.
Let $u$ be a non null temperated distribution.
If $H_{u,\alpha }$ is bounded from $L^{p_{1},q_{1}}\times L^{p_{2},q_{2}}$
into $L^{p_{3},q_{3}}$  with norm $\| H_{u,\al }\| $ then
$0\leq \frac{1}{p_{1}}+\frac{1}{p_{2}}-\frac{1}{p_{3}}\leq 1$.

In this
case, if $G(x)=e^{-\pi x^2}$ and
$\frac{1}{p}=\frac{1}{p_{1}}+\frac{1}{p_{2}}-\frac{1}{p_{3}}$ we have that
$\hat{u}*D_{\lambda }^{p'}G$ is a uniformly bounded family of functions with
$$
\sup_{\lambda >0}\| \hat{u}*D_{\lambda }^{p'}G\|_{\infty }
\leq C\| H_{u,\al }\|
$$
where $C$ is a constant that depends only of $\al $, $p_{i}$ and $q_{i}$,
$i=1,2,3$.
\end{prop}
\begin{nota} When $\frac{1}{p_{1}}+\frac{1}{p_{2}}=\frac{1}{p_{3}}$
the thesis says that $\hat{u}$ is a bounded function with
$
\| \hat{u}\|_{\infty }\leq C\| H_{u,\al }\|
$ which is a known fact for linear multipliers (see \cite{larsen}).
\end{nota}
\Dem
Let $\omega \in
\bbbr $, $\lambda >0$, $\alpha \in \bbbr \backslash \{ 0,1\} $ and
define $\lambda'=(1+|\alpha |)^{-1}\lambda^{2}$.
Let $f(t)=e^{2\pi i\omega t}e^{-\lambda' \pi t^2}$ and
$g(t)
=e^{-\frac{\lambda' }{|\alpha |}\pi t^2}$. An easy computation shows that
%$$
%f(x-t)g(x-\al t)
%=e^{2\pi i\omega (x-t)}e^{-\lambda' \pi (x-t)^2}
%e^{-\frac{\lambda' }{|\alpha |}\pi (x-\alpha t)^2}
%$$
%$$
%=e^{2\pi i\omega x}e^{-2\pi i\omega t}
%e^{-\lambda' \pi x^2}e^{2\lambda' \pi xt}e^{-\lambda' \pi t^2}
%e^{-\frac{\lambda' }{|\alpha |}\pi x^2}
%e^{2sign{(\alpha )}\lambda' \pi xt}e^{-|\alpha |\lambda' \pi t^2}
%$$
for $\alpha <0$ we have
$
f(x-t)g(x-\al t)
%=e^{2\pi i\omega x}
%e^{-\lambda' \pi x^2}e^{-\frac{\lambda' }{|\alpha |}\pi x^2}
%e^{-2\pi i\omega t}e^{-\lambda' \pi t^2}
%e^{-|\alpha |\lambda' \pi t^2}
=f(x)g(x)f(-t)g(-\alpha t)
$.
Thus
$$
H_{u,\al }(f,g)(x)=f(x)g(x)H_{u,\alpha }(f,g)(0)
$$
which says that the BHT of these gaussian functions is
the product of both functions times a constant. Since
$$
\big|H_{u,\alpha }(f,g)(0)\big|\| fg\|_{p_{3},q_{3}}
=\| H_{u,\alpha }(f,g)\|_{p_{3},q_{3}}
\leq \| H_{u,\al }\| \| f\|_{p_{1},q_{1}}\| g\|_{p_{2},q_{2}}
$$
we just need to compute norms
in order to get the desired  condition:
$$
\| f\|_{p_{1},q_{1}}=\| M_{\omega }D_{{\lambda'}^{-1/2}}G\|_{p_{1},q_{1}}
={\lambda'}^{-\frac{1}{2p_{1}}}\| G\|_{p_{1},q_{1}}
$$
$$
\| g\|_{p_{2},q_{2}}
=\| D_{(\frac{\lambda'}{|\alpha |})^{-1/2}}G\|_{p_{2},q_{2}}
={\lambda'}^{-\frac{1}{2p_{2}}}|\alpha |^{\frac{1}{2p_{2}}}
\| G\|_{p_{2},q_{2}}
$$
$$
\| fg\|_{p_{3},q_{3}}=\| M_{\omega }
D_{{\lambda'}^{-\frac{1}{2}}(1+\frac{1}{|\al
|})^{-\frac{1}{2}}}G\|_{p_{3},q_{3}}  ={\lambda'}^{-\frac{1}{2p_{3}}} \Big(
1+\frac{1}{|\alpha |}
\Big)^{-\frac{1}{2p_{3}}}\| G\|_{p_{3},q_{3}}
$$
with
$$
\| G\|_{p_{i},q_{i}}
=\bigg( \frac{q_{i}}{2p_{i}}\bigg)^{\frac{1}{q_{i}}}
\Gamma \bigg( \frac{q_{i}}{2p_{i}}\bigg)^{\frac{1}{q_{i}}}
\bigg( \frac{4}{q_{i}\pi }\bigg)^{\frac{1}{2p_{i}}}
$$
where $\Gamma $ denotes the function Gamma of Euler (see remark
\ref{lorentznorm} below). So
$$
\big| H_{u,\alpha }(f,g)(0)\big|
\leq \| H_{u,\al }\| \frac{\| G\|_{p_{1},q_{1}}\| G\|_{p_{2},q_{2}}}{\|
G\|_{p_{3},q_{3}}} |\alpha |^{{1}\over{2p_{2}}}
\bigg( 1+\frac{1}{|\alpha |}\bigg)^{\frac{1}{2p_{3}}}
{\lambda'}^{-\frac{1}{2}\big(
\frac{1}{p_{1}}+\frac{1}{p_{2}}-\frac{1}{p_{3}}\big) }
$$
$$
=\| H_{u,\al }\|\frac{\| G\|_{p_{1},q_{1}}\| G\|_{p_{2},q_{2}}}{\|
G\|_{p_{3},q_{3}}}
|\al |^{\12 \big( \frac{1}{p_{2}}-\frac{1}{p_{3}}\big) }
(1+|\alpha |)^{\12 \big( \frac{1}{p_{1}}+\frac{1}{p_{2}}\big) }
%|\alpha |^{{1}\over{2p_{2}}}\bigg( 1+\frac{1}{|\alpha
%|}\bigg)^{\frac{1}{2p_{3}}} (1+|\al |)^{\frac{1}{2p}}
\lambda^{-\frac{1}{p}}
=C \lambda^{-\frac{1}{p}}
$$
for all $\lambda >0$ and $\omega \in \bbbr $.
%If $q^{-1}=q_{1}^{-1}+q_{2}^{-1}-q_{3}^{-1}$ the value of the constant is
%$$
%C=\frac{2^{\frac{1}{p}}}
%{2^{\frac{1}{q}}}
%\frac{p_{3}^{\frac{1}{q_{3}}}}
%{p_{1}^{\frac{1}{q_{1}}}p_{2}^{\frac{1}{q_{2}}}}
%\frac{q_{1}^{\frac{1}{q_{1}}-\frac{1}{2p_{1}}}
%q_{2}^{\frac{1}{q_{2}}-\frac{1}{2p_{2}}}}
%{q_{3}^{\frac{1}{q_{3}}-\frac{1}{2p_{3}}}}
%\Gamma \bigg( \frac{q_{1}}{2p_{1}}\bigg)^{\frac{1}{q_{1}}}
%\Gamma \bigg( \frac{q_{2}}{2p_{2}}\bigg)^{\frac{1}{q_{2}}}
%\Gamma \bigg( \frac{q_{3}}{2p_{3}}\bigg)^{-\frac{1}{q_{3}}}
%$$
%$$
%|\al |^{\12 \big( \frac{1}{p_{2}}-\frac{1}{p_{3}}\big) }
%(1+|\alpha |)^{\12 \big( \frac{1}{p_{1}}+\frac{1}{p_{2}}\big) }
%\| H_{u,\al }\|
%$$

Now we work a little bit the expression $H_{u,\alpha }(f,g)(0)$. Since
$$
f(-t)g(-\alpha t)=e^{-2\pi i\omega t}e^{-(1+|\alpha |)\lambda' \pi t^2}
=e^{-2\pi i\omega t}e^{-\lambda^2 \pi t^2}
=M_{-\omega }D_{\lambda^{-1}}G(t)
$$
we have, using the fact that $\widehat{G}=G$ and $D_{-1}G=G$, that
\begin{equation}\label{conv}
H_{u,\alpha }(f,g)(0)=(u,M_{-\omega }D_{\lambda^{-1}}G)
%=(\hat{u},(M_{-\omega }D_{\lambda^{-1}}G)\checke{})
%$$
=(\hat{u},T_{\omega }D_{\lambda }^{1}G)
%=(\hat{u},T_{\omega }D_{-1}D_{\lambda }^{1}G)
=(\hat{u}*D_{\lambda }^{1}G)(\omega )
\end{equation}
and we can rewrite the previous result for all $\lambda >0$ and $\omega \in \bbbr $ as
$$
|(\hat{u}*D_{\lambda }^{1}G)(\omega )|
\leq C\lambda^{-\frac{1}{p}}
$$

a) If $\frac{1}{p}<0$
we prove that $u\equiv 0$ by showing that the family of functions
$
m_{\lambda }(\omega )=(\hat{u}*D_{\lambda }^{1}G)(\omega )
$
converge pointwise to zero and distributionally to $\hat{u}$ when $\lambda $
tends to zero.

On one side, we have $m_{\lambda }$ are bounded functions
(and so locally integrable) with
$
\| m_{\lambda }\|_{\infty }
\leq C\lambda^{-\frac{1}{p}}\leq C
$
for $\lambda <1$ and
$
\lim_{\lambda \rightarrow 0}m_{\lambda }(\omega )=0
$
for all $\omega \in \bbbr $.

On the other side,
since $(D_{\lambda }^{1}G)_{\lambda >0}$ is an approximate identity
we have proven in lemma \ref{aproxid} that
$\{ D_{\lambda }^{1}G*\varphi )\}_{\lambda >0}$
converges to $\varphi $ in the topology $\T_{\S }$. Thus,
by continuity of $\hat{u}$ we have for all $\varphi \in \S $
$$
\lim_{\lambda \rightarrow 0}(u_{m_{\lambda }},\varphi )
=\lim_{\lambda \rightarrow 0}(\hat{u}*D_{\lambda }^{1}G,\varphi )
=\lim_{\lambda \rightarrow 0}(\hat{u},D_{\lambda }^{1}G*\varphi )
=(\hat{u},\varphi )
$$

With both facts and dominated convergence theorem of Lebesgue we have
%for all $\varphi \in \S $
$$
(\hat{u},\varphi )=\lim_{\lambda \rightarrow 0}(u_{m_{\lambda }},\varphi )
=\lim_{\lambda \rightarrow 0}\int_{\re }m_{\lambda }(\omega )\varphi (\omega )
d\omega =0
$$

b) If $\frac{1}{p}=0$
we still know that
$
m_{\lambda }
$
define a family of bounded functions with $\| m_{\lambda }\|_{\infty }\leq
C$ for all $\lambda >0$ that converge distributionally to
$\hat{u}$ when $\lambda $ tends to zero. We use this fact to show that
$\hat{u}$ must be a bounded function and that, actually, the convergence is
also pointwise.
From above,
$$
|(\hat{u},\varphi )|
%=\lim_{\lambda \rightarrow 0}|(u_{m_{\lambda }},\varphi )|
=\lim_{\lambda \rightarrow 0}
\Big| \int_{\re } m_{\lambda }(\omega )\varphi (\omega )d\omega \Big|
\leq \overline{\lim_{\lambda \rightarrow 0}}
\| m_{\lambda }\|_{\infty }\| \varphi \|_{1}
\leq C\| \varphi \|_{1}
$$
for all $\varphi \in \S $ and thus $\hat{u}$ is a distribution associated to a bounded function.
%with bound less or equal to $C$.
%Since $\hat{u}$ is a bounded function and
%$
%m_{\lambda }(\omega )
%5=(\hat{u}*D_{\lambda }^{1}G)(\omega )
%$
%we know,
Moreover, by property of approximate identity, we have that
$$
\lim_{\lambda \rightarrow 0}m_{\lambda }(\omega )
=\lim_{\lambda \rightarrow 0}(\hat{u}*D_{\lambda }^{1}G)(\omega )
=\hat{u}(\omega )
$$
almost everywhere (in all Lebesgue points of $\hat{u}$).

c) If $0<\frac{1}{p}\leq 1$
our  condition says that
$
|(\hat{u}*D_{\lambda }^{p'}G)(\omega )|\leq C
$
for all $\lambda >0$ and $\omega \in \bbbr $ which is the main statement of
the proposition.

We still have that
$
m_{\lambda }=\hat{u}*D_{\lambda }^{1}G
$
define a family of bounded functions that converge distributionally to
$\hat{u}$ and satisfies $\| m_{\lambda }\|_{\infty }\leq
C\lambda^{-\frac{1}{p}}$ for all $\lambda >0$.

d) If $1<\frac{1}{p}$ we prove directly that $u\equiv 0$. The previous condition can be written as
$
|(\hat{u}*D_{\lambda }G)(\omega )|\leq C\lambda^{\frac{1}{p'}}
$
with $p'<0$. Moreover, since $H_{u,\alpha }$ is bounded and translation invariant by property (\ref{invtrasl}),
we have that
$H_{T_{y}u,\alpha }$ is also a bounded operator with the same constant and thus it satisfies
$
|(\widehat{T_{y}u}*D_{\lambda }G)(\omega )|\leq C\lambda^{\frac{1}{p'}}
$
for every $y,\omega \in \bbbr $. With this we can write
$$
\lim_{\lambda \rightarrow 0}|(u,T_{y}D_{\lambda }^{1}G)|
%=\lim_{\lambda \rightarrow 0}|(T_{-y}u,D_{\lambda }^{1}G)|
=\lim_{\lambda \rightarrow 0}|(\widehat{T_{-y}u},D_{\lambda^{-1} }G)|
$$
$$
=\lim_{\lambda \rightarrow 0}|(\widehat{T_{-y}u}*D_{\lambda^{-1}}G)(0)|
\leq \lim_{\lambda \rightarrow 0}C\lambda^{-\frac{1}{p'}}=0
$$

Thus for every $\varphi \in \S $ we have by the dominated convergence theorem
$$
(u,\varphi )=\lim_{\lambda \rightarrow 0}(u,\varphi *D_{\lambda }^{1}G)
=\lim_{\lambda \rightarrow 0}\int_{\re }\varphi (y)
(u,T_{y}D_{\lambda }^{1}G)dy=0
$$

Now we see the case when $\alpha >0$. If $\al >1$ and $p_{3}\geq 1$ the duality formula
(\ref{dua}) with $1-\al <0$ let us to apply the former result to
$H_{D_{-1}u,1-\alpha }$   in the following way: if $f,g,h$ are some properly
chosen gaussian functions then
$$
\big\langle h,H_{u,\alpha }(f,g)\big\rangle
=\big\langle H_{D_{-1}u,1-\alpha }(h,g),f\big\rangle
=H_{D_{-1}u,1-\alpha }(h,g)(0)\big\langle hg,f \big\rangle
$$
which, if we claim the operator bounded, implies
$$
|H_{D_{-1}u,1-\alpha }(h,g)(0)|\leq \| H_{u,\alpha }\|
\frac{\| f\|_{p_{1},q_{1}}\| g\|_{p_{2},q_{2}}\| h\|_{p_{3}',q_{3}'}}{\big|
\langle hg,f\rangle \big|}
=C\lambda^{-\frac{1}{p}}
$$
Thus by (\ref{conv}) and using $D_{-1}(f*g)=D_{-1}f*D_{-1}g$,
$D_{-1}\hat{u}=\widehat{D_{-1}u}$ we have
$$
|(\hat{u}*D_{\lambda }^{1}G)(-\omega )|
=|(\widehat{D_{-1}u}*D_{\lambda }^{1}G)(\omega )|
=|H_{D_{-1}u,1-\alpha }(h,g)(0)|\leq C\lambda^{-\frac{1}{p}}
$$
From here the same ideas lead to the same conclusion.

Finally, if $0<\al <1$ and $p_{3}\geq 1$, the commutativity formula
(\ref{conmu}) with $\al^{-1}>1$ and the duality formula (\ref{dua}) with
$1-\alpha^{-1}<0$
let us apply the same ideas to
$H_{D_{-1}D_{\al }^{1}u,1-\alpha^{-1}}$ to get the same conclusion:
$$
\big\langle h,H_{u,\alpha }(f,g)\big\rangle
=\big\langle h,H_{D_{\alpha }^{1}u,\alpha^{-1} }(g,f)\big\rangle
$$
$$
=\big\langle H_{D_{-1}D_{\alpha }^{1}u,1-\alpha^{-1}}(h,f),g\big\rangle
=H_{D_{-1}D_{\alpha }^{1}u,1-\alpha^{-1}}(h,f)(0)
\big\langle hf,g\big\rangle
$$
which implies
$$
|H_{D_{-1}D_{\alpha }^{1}u,1-\alpha^{-1}}(h,f)(0)|\leq \| H_{u,\alpha }\|
\frac{\| f\|_{p_{1},q_{1}}\| g\|_{p_{2},q_{2}}\| h\|_{p_{3}',q_{3}'}}{\big|
\langle hf,g\rangle \big|}
=C\lambda^{-\frac{1}{p}}
$$
Now, using that $D_{\al }(f*g)=D_{\al }^{q}f*D_{\al }^{q'}g$ we get
by (\ref{conv})
$$
|(\hat{u}*D_{\al \lambda }^{1}G)(-\al \omega )|
=|D_{\alpha^{-1}}(\hat{u}*D_{\al }^{1}D_{\lambda }^{1}G)(-\omega )|
=|(D_{\al^{-1} }\hat{u}*D_{\lambda }^{1}G)(-\omega )|
$$
$$
%=|(\widehat{D_{\alpha }^{1}u}*D_{\lambda }^{1}G)(-\omega ))|
=|(\widehat{D_{-1}D_{\alpha }^{1}u}*D_{\lambda }^{1}G)(\omega )|
=|H_{D_{-1}D_{\alpha }^{1}u,1-\alpha^{-1}}(h,f)(0)|\leq
C\lambda^{-\frac{1}{p}}
$$
and we finish with the same ideas as before.

%\begin{nota}
%We have actually proved that if $u$ is a temperated distribution
%such that $H_{u,\al }$ is a bounded operator for a particular value of
%$\al $ then $\hat{u}*D_{\lambda }^{p}G$ is a uniformly bounded family
%of functions. This is because from (\ref{conv}) we have that
%$
%(\hat{u}*D_{\lambda }^{1}G)(\omega )=H_{u,\alpha }(f,g)(0)
%$
%for
%$f(t)=e^{2\pi i\omega t}e^{-\lambda' \pi t^2}$ and
%$g(t)=e^{-\frac{\lambda' }{|\alpha |}\pi t^2}$ where
%$\lambda' =(1+|\alpha |)^{-1}\lambda^{2}$.

%\end{nota}
\begin{nota}\label{lorentznorm}
Since
$G$ is even and non-increasing in $[0,\infty )$, we know that
$G^{*}=D_{2}G$ and so we can compute
$
\| G\|_{p_{1},q_{1}}^{q_{1}}
$
as follows
$$
\frac{q_{1}}{p_{1}}\int_{0}^{\infty }
t^{\frac{q_{1}}{p_{1}}}e^{-\frac{1}{4}q_{1}\pi t^2}\frac{dt}{t}
%=\frac{q_{1}}{2p_{1}}\int_{0}^{\infty }
%t^{\frac{q_{1}}{2p_{1}}}e^{-\frac{1}{4}q_{1}\pi t}\frac{dt}{t}
%$$
%$$
=\frac{q_{1}}{2p_{1}}\bigg( \frac{4}{q_{1}\pi
}\bigg)^{\frac{q_{1}}{2p_{1}}}
\int_{0}^{\infty }t^{\frac{q_{1}}{2p_{1}}}e^{-t}\frac{dt}{t}
=\frac{q_{1}}{2p_{1}}\bigg( \frac{4}{q_{1}\pi
}\bigg)^{\frac{q_{1}}{2p_{1}}}
\Gamma \bigg( \frac{q_{1}}{2p_{1}}\bigg)
$$
\end{nota}

\subsection{Bilinear Young inequality.}
The next result is the generalization
of Young inequality to our bilinear non-convolution operators. We pay now
special attention to the dependence of the constants from the
parameter $\alpha $. In order to deal with a more general and
symmetric operator, we change a little bit its definition. For the next
proposition we call BHT to
$$
H_{K,\alpha ,\beta }(f,g)(x)=\int_{\re}f(x-\alpha t)g(x-\beta t)K(t)dt
$$
defined for all $\alpha, \beta , x\in \bbbr $ and $f,g\in \S $.

\begin{prop}{\bf (Bilinear Young inequality).}\label{youngbilineal}
Let $p_{0}\geq 1$.
If $K\in L^{p_{0}}$
then $H_{K,\alpha ,\beta }$ is a bounded operator from
$L^{p_{1}}\times L^{p_{2}}$ to  $L^{p_{3}}$
with $p_{i}\geq 1$ for $i=1,2,3$ and $p_{1}^{-1}+p_{2}^{-1}+p_{0}^{-1}=1+{p_{3}}^{-1}$,
%\begin{equation}\label{expo}
%\frac{1}{p_{1}}+\frac{1}{p_{2}}+\frac{1}{p_{0}}=1+\frac{1}{p_{3}}
%\end{equation}
and all $\alpha ,\beta \in \bbbr \backslash \{ 0\} $ such that $\alpha \neq
\beta
$. Moreover,
$$
\| H_{K,\alpha ,\beta }(f,g)\|_{p_{3}}
\leq C_{\alpha ,\beta ,p_{0},p_{1},p_{2}}\| K\|_{p_{0}} \|f\|_{p_{1}}
\|g\|_{p_{2}}
$$
\end{prop}
\begin{nota} Notice that
$
p_{1}^{-1}+p_{2}^{-1}-p_{3}^{-1}={p_{0}'}^{-1}\in [0,1]
$
as proposition \ref{gaussian} says it must be.
See also that this condition can be rewritten as
$
p_{1}^{-1}+p_{2}^{-1}+{p_{3}'}^{-1}=1+{p_{0}'}^{-1}
$
and so one can think the point $(p_{1}^{-1},p_{2}^{-1},{p_{3}'}^{-1})\in
\bbbr^{3}$ belongs to the plane
$x+y+z=1+{p_{0}'}^{-1}$ with $1+{p_{0}'}^{-1}\in [1,2]$.
\end{nota}
\Dem Let $p\geq 1$, $f,g,h,K\in \S $ and
$
I=\Big| \int_{\re }h(x)\int_{\re }f(x-\al t)g(x-\beta t)K(t)dtdx\Big|
$.
We denote here $f_{a,b}(x,t)=f(ax+bt)$.
By H\"older inequality and some changes of variables\\
$
I
\leq \| f_{1,-\alpha }g_{1,-\beta }\|_{L^{p}(\bbbr^{2})}
\| K_{0,1}h_{1,0}\|_{L^{p'}(\bbbr^{2})}
=|\alpha -\beta |^{-\frac{1}{p}}
\| f\|_{p}\| g\|_{p}\| K\|_{p'}\| h\|_{p'}
$ \ecu i.e.
\begin{equation}\label{ab}
\| H_{K,\alpha ,\beta }(f,g)\|_{p}
\leq |\alpha -\beta |^{-\frac{1}{p}}\| f\|_{p}\| g\|_{p}\|
K\|_{p'}
\end{equation}
$
I\leq \| K_{0,1}g_{1,-\beta }\|_{L^{p}(\bbbr^{2})}
\| f_{1,-\alpha }h_{1,0}\|_{L^{p'}(\bbbr^{2})}
=|\alpha |^{-\frac{1}{p'}}\| f\|_{p'}\| g\|_{p}\| K\|_{p}\| h\|_{p'}
$\ecu i.e.
\begin{equation}\label{a}
\| H_{K,\alpha ,\beta }(f,g)\|_{p}
\leq |\alpha |^{-\frac{1}{p'}} \| f\|_{p'}\| g\|_{p}\| K\|_{p}
\end{equation}
$
I\leq \| f_{1,-\alpha }K_{0,1}\|_{L^{p}(\bbbr^{2})}
\| g_{1,-\beta }h_{1,0}\|_{L^{p'}(\bbbr^{2})}
=|\beta |^{-\frac{1}{p'}}\| f\|_{p}\| g\|_{p'}\| K\|_{p}\| h\|_{p'}
$\ecu i.e.
\begin{equation}\label{b}
\| H_{K,\alpha ,\beta }(f,g)\|_{p}
\leq |\beta |^{-\frac{1}{p'}} \| f\|_{p}\| g\|_{p'}\| K\|_{p}
\end{equation}

We associate each bound of the operator from $L^{p_{1}}\times L^{p_{2}}$
to $L^{p_{3}}$ to the point
$(p_{1}^{-1},p_{2}^{-1},{p_{3}'}^{-1})\in \bbbr^{3}$ in the
plane $x+y+z=1+p^{-1}$. In this way and taking the values $p=1$ and
$p=\infty
$ in each of the three  previous inequalities we consider
the extremal points $(1,1,0)$, $(0,0,1)$ (from the first one), $(0,1,0)$,
$(1,0,1)$ (from the second), $(1,0,0)$ and $(0,1,1)$ (from the third).
In this way, by using trilinear interpolation between two spaces iteratively we
get the bounds on the surface of
the convex hull of the previous six points, that is, on the surface of the
octahedron drawn in the following diagram

\begin{center}
\begin{picture}(160,160)
%\put(0,0){\framebox(170,170)}
\put(15,5){\includegraphics[scale=.7]{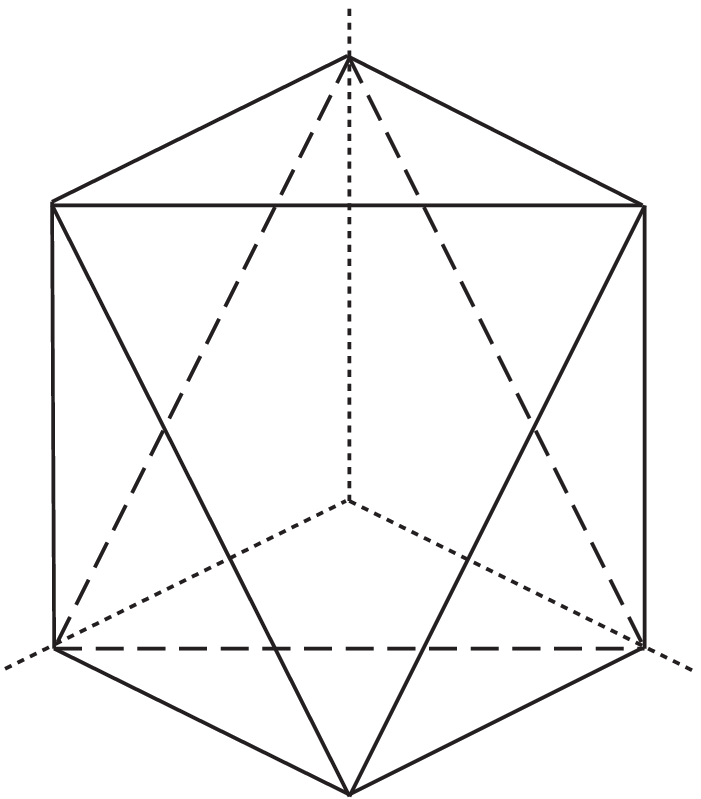}}
\put(5,125){\tiny $|\alpha |^{-1}$}
\put(150,125){\tiny $|\beta |^{-1}$}
\put(70,0){\tiny $|\alpha -\beta |^{-1}$}
\put(22,30){\tiny $1$}
\put(145,30){\tiny $1$}
\put(80,155){\tiny $1$}
\put(70,20){\tiny $|\alpha -\beta |^{-\frac{1}{p_{0}'}}$}
\put(130,50){\tiny $\leftarrow
\!\!\!\!-\!\!\!-\!\!\!-\!\!\!-\!\!\!-\!\!\!-\!\!\!-
|\beta |^{-\frac{1}{p_{3}'}} |\alpha -\beta |^{-\frac{1}{p_{1}}}$}
\put(130,80){\tiny $\leftarrow - - - |\beta |^{-\frac{1}{p_{0}'}}$}
\put(-60,50){\tiny $|\alpha |^{-\frac{1}{p_{3}'}}
|\alpha -\beta |^{-\frac{1}{p_{2}}}
-\!\!\!\!-\!\!\!-\!\!\!-\!\!\!-\!\!\!-\!\!\!-\!\!\!-\!\!\!\rightarrow$}
\put(-20,80){\tiny $|\alpha |^{-\frac{1}{p_{0}'}} - - - \rightarrow$}
\put(62,130){\tiny $|\alpha |^{-\frac{1}{p_{1}}}
|\beta |^{-\frac{1}{p_{2}}}$}
\put(43,95){\tiny $|\alpha |^{-\frac{1}{p_{2}'}}
|\beta |^{-\frac{1}{p_{1}'}} |\alpha -\beta |^{-\frac{1}{p_{3}}}$}
\put(15,5){\includegraphics[scale=.7]{diagrama.eps}}
\end{picture}
\end{center}
where we write the constants of boundedness in each vertex and each
face. We show how to get one of them: from (\ref{a}) and (\ref{b}) we know
$
\| H_{K,\alpha ,\beta }(f,g)\|_{\infty }
\leq |\alpha |^{-1} \| f\|_{1}\| g\|_{\infty }\| K\|_{\infty }
$, $\| H_{K,\alpha ,\beta }(f,g)\|_{\infty }
\leq |\beta |^{-1} \| f\|_{\infty }\| g\|_{1}\| K\|_{\infty }
$
and so we have
$
\| H_{K,\alpha ,\beta }(f,g)\|_{\infty  }
\leq |\alpha |^{-\frac{1}{p}}|\beta |^{-\frac{1}{p'}}
\| f\|_{p}\| g\|_{p'}\| K\|_{\infty }
$.
In the same way, from (\ref{a}) and (\ref{ab})
$
\| H_{K,\alpha ,\beta }(f,g)\|_{\infty }
\leq |\alpha |^{-1} \| f\|_{1}\| g\|_{\infty }\| K\|_{\infty }
$,
$
\| H_{K,\alpha ,\beta }(f,g)\|_{\infty }
\leq \| f\|_{\infty }\| g\|_{\infty }\| K\|_{1}
$
we get
$
\| H_{K,\alpha ,\beta }(f,g)\|_{\infty  }
\leq |\alpha |^{-\frac{1}{p}}
\| f\|_{p}\| g\|_{\infty }\| K\|_{p'}
$.

Interpolating both cases we get
$
\| H_{K,\alpha ,\beta }(f,g)\|_{\infty  }
\leq |\alpha |^{-\frac{1}{p}}|\beta |^{-\frac{1}{q_1}}
\| f\|_{p}\| g\|_{q_1}\| K\|_{q_2}
$
with $q_1^{-1}+q_2^{-1}={p'}^{-1}$.
Using again (\ref{ab})
$
\| H_{K,\alpha ,\beta }(f,g)\|_{1}
\leq |\alpha -\beta |^{-1}
\| f\|_{1}\| g\|_{1}\| K\|_{\infty }
$,
we finally have
$$
\| H_{K,\alpha ,\beta }(f,g)\|_{p_{3}}
\leq |\alpha |^{-\frac{1}{p_1}}|\beta |^{-\frac{1}{p_2}}
|\alpha -\beta |^{-\frac{1}{p_3}}
\| f\|_{p_{1}}\| g\|_{p_{2}}\| K\|_{p_{0}}
$$
where $p_{3}^{-1}=\theta $, $p_{1}^{-1}=(1-\theta )p^{-1}+\theta $,
$p_{2}^{-1}=(1-\theta ){q_1}^{-1}+\theta $ and
$p_{0}^{-1}=(1-\theta ){q_2}^{-1}$,
which is the stated result since $p_{1}^{-1}+p_{2}^{-1}+{p_{0}}^{-1}=1+p_3^{-1}$.

Now in order to get bounds in the interior of the octahedron we
use interpolation between six spaces. In this way, each point
$p=(p_{1}^{-1},p_{2}^{-1},{p_{3}'}^{-1})$ can be written as the convex linear
combination of the six vertex in the following way
$$
p=(\lambda_{2}+p_{3}^{-1}-p_{2}^{-1})(1,0,0)
+(\lambda_{1}+p_{3}^{-1}-p_{1}^{-1})(0,1,0)
+({p_{3}'}^{-1}-\lambda_{1}-\lambda_{2})(0,0,1)
$$
$$
+\lambda_{1}(1,0,1)
+\lambda_{2}(0,1,1)
+(p_{1}^{-1}+p_{2}^{-1}-p_{3}^{-1}-\lambda_{1}-\lambda_{2})(1,1,0)
$$
for every $\lambda_{1},\lambda_{2}\in [0,1]$ such that
$\max(p_{1}^{-1}-p_{3}^{-1},0)\leq \lambda_{1}$,
$\max(p_{2}^{-1}-p_{3}^{-1},0)\leq \lambda_{2}$ and
$\lambda_{1}+\lambda_{2}\leq \min({p_{3}'}^{-1},{p_{0}'}^{-1})$.
We denote by $D$ such non empty triangle
(notice that $p_{i}^{-1}-p_{3}^{-1}\leq
p_{1}^{-1}+p_{2}^{-1}-p_{3}^{-1} ={p_{0}'}^{-1}\leq1$ and
$\max(p_{1}^{-1}-p_{3}^{-1},0)+
\max(p_{2}^{-1}-p_{3}^{-1},0)\leq \min({p_{3}'}^{-1},{p_{0}'}^{-1})$).
Also notice that this decomposition implies this other one
for $\tilde{p}=(p_{1}^{-1},p_{2}^{-1},p_{3}^{-1})$
$$
\tilde{p}
=(\lambda_{2}+p_{3}^{-1}-p_{2}^{-1})(1,0,1)
+(\lambda_{1}+p_{3}^{-1}-p_{1}^{-1})(0,1,1)
+(1-{p_{3}'}^{-1}+\lambda_{1}+\lambda_{2})(0,0,0)
$$
$$
+\lambda_{1}(1,0,0)
+\lambda_{2}(0,1,0)
+(p_{1}^{-1}+p_{2}^{-1}-p_{3}^{-1}-\lambda_{1}-\lambda_{2})(1,1,1)
$$
in order to interpolate. So, using theorem \ref{interpolate} we get
$$
\| H_{K,\alpha ,\beta }(f,g)\|_{p_{3}}
\leq |\alpha |^{-\lambda_{1}}|\beta |^{-\lambda_{2}}
|\alpha -\beta |^{-(\frac{1}{p_{0}'}-\lambda_{1}-\lambda_{2})}
\| f\|_{p_{1}}\| g\|_{p_{2}}\| K\|_{p_{0}}
$$
for every $\lambda_{1},\lambda_{2}\in D$ and we want now to minimize.
Since $D$ is a convex domain and
$F(x,y)=(|\alpha ||\alpha -\beta |^{-1})^{-x}(|\beta ||\alpha -\beta
|^{-1})^{-y}$ is a convex function in $D$, the minimal costant is attained
in one of the three vertex of the triangle:
$$\hspace{-1.2cm}
(\max(p_{1}^{-1}-p_{3}^{-1},0),\max(p_{2}^{-1}-p_{3}^{-1},0))
$$
$$
(\max(p_{1}^{-1}-p_{3}^{-1},0),
\min({p_{1}'}^{-1},p_{2}^{-1},{p_{3}'}^{-1},{p_{0}'}^{-1}))
$$
$$
(\min(p_{1}^{-1},{p_{2}'}^{-1},{p_{3}'}^{-1},{p_{0}'}^{-1}),
\max(p_{2}^{-1}-p_{3}^{-1},0))
$$
that is
$$
\| H_{K,\alpha ,\beta }(f,g)\|_{p_{3}}
\leq C_{\alpha ,\beta ,p_{1},p_{2},p_{0}}
\| f\|_{p_{1}}\| g\|_{p_{2}}\| K\|_{p_{0}}
$$
where $C_{\alpha ,\beta ,p_{1},p_{2},p_{0}}$ is the minimum of the
three quantities:
$$
\hspace{-.4cm}
|\alpha |^{-\max(p_{1}^{-1}-p_{3}^{-1},0)}
|\beta |^{-\max(p_{2}^{-1}-p_{3}^{-1},0)}
|\alpha -\beta |^{-\min(p_{1}^{-1},p_{2}^{-1},p_{3}^{-1},{p_{0}'}^{-1})}
$$
$$
|\alpha |^{-\max(p_{1}^{-1}-p_{3}^{-1},0)}
|\beta |^{-\min({p_{1}'}^{-1},p_{2}^{-1},{p_{3}'}^{-1},{p_{0}'}^{-1})}
|\alpha -\beta |^{-\max(p_{2}^{-1}-{p_{1}'}^{-1},0)}
$$
$$
|\alpha |^{-\min(p_{1}^{-1},{p_{2}'}^{-1},{p_{3}'}^{-1},{p_{0}'}^{-1})}
|\beta |^{-\max(p_{2}^{-1}-p_{3}^{-1},0)}
|\alpha -\beta |^{-\max(0,p_{1}^{-1}-{p_{2}'}^{-1})}
$$
which, on the surface of the octahedron, are the same bounds we
already had (in fact, the three bounds coincide in each face).
\cqd

\subsection{The third condition} The last result gives a sufficient condition
of boundedness for bilinear multipliers.
It gives a condition over the symbol of the operator instead of over
the kernel.
\begin{prop}
Let $m\in L^{q}(\bbbr^{2})$ with $1\leq q\leq 4$. Then $m$ is
$(p_{1},p_{2},p_{3})$-multiplier for all exponents such that
$1\leq p_{1},p_{2},p_{3}'\leq \min (2,q)$,
$q\notin \{ p_{1},p_{2},p_{3}'\} $ and
$p_{1}^{-1}+p_{2}^{-1}+{p_{3}'}^{-1}=1+2q^{-1}$. Moreover, $\|
m\|_{\MB_{p_{1},p_{2},p_{3}}}\leq \| m\|_{q}$.
\end{prop}

\Dem
By duality it is enough to prove that
for every $f,g,h\in \S $
$$
I=\Big| \int_{\re^{2}}\hat{f}(\xi )\hat{g}(\eta )m(\xi ,\eta )
\hat{h}(-\xi -\eta )d\xi d\eta \Big|
\leq C_{m}\| f\|_{p_{1}}\| g\|_{p_{2}}\| h\|_{p_{3}'}
$$

If $q=1$ then
%$p_{1}=p_{2}=p_{3}'=1$ and thus
$
%\Big| \int_{\re^{2}}\hat{f}(\xi )\hat{g}(\eta )m(\xi ,\eta )
%\hat{h}(-\xi -\eta )d\xi d\eta \Big|\leq
I\leq \| m\|_{1}\| \hat{f}\|_{\infty } \| \hat{g}\|_{\infty }\|
\hat{h}\|_{\infty }
\leq \| m\|_{1}\| f\|_{1}\| g\|_{1}\| h\|_{1}
$.

If $q>1$, we define $\tilde{p}=(\tilde{p}_{1},\tilde{p}_{2},\tilde{p}_{3}')$
by
$$
\tilde{p}_{1}=\frac{p_{1}(q-1)}{q-p_{1}}
\e \e
\tilde{p}_{2}=\frac{p_{2}(q-1)}{q-p_{2}}
\e \e
\tilde{p}'_{3}=\frac{p_{3}'(q-1)}{q-p_{3}'}
$$
which satisfy:
$$
1\leq \tilde{p}_{1},\tilde{p}_{2},\tilde{p}_{3}'\leq \infty ,
\hskip 15pt
\tilde{p}_{i}'=\frac{p_{i}'}{q'}
\hskip 10pt
i=1,2 ,
\hskip 15pt
\tilde{p}_{3}=\frac{p_{3}}{q'},
\hskip 15pt
\frac{1}{\tilde{p}'_{1}}+\frac{1}{\tilde{p}'_{2}}+\frac{1}{\tilde{p}_{3}}=2
$$

Then, by Hölder, Young and Hausdorff-Young inequalities we have
%$$
%\Big| \int_{\re^{2}}\hat{f}(\xi )\hat{g}(\eta )m(\xi ,\eta )
%\hat{h}(-\xi -\eta )d\xi d\eta \Big|
%$$
$$
I\leq \| m\|_{q}
\Big( \int_{\re^{2}}|\hat{f}(\xi )|^{q'}|\hat{g}(\eta )|^{q'}
|\hat{h}(-\xi -\eta )|^{q'}d\xi d\eta \Big)^{\frac{1}{q'}}
%$$
%$$
%=\| m\|_{q}
%\Big( \int_{\re }|\hat{f}(\xi )|^{q'}
%(|\hat{g}|^{q'}*|\hat{h}|^{q'})(-\xi )d\xi\Big)^{\frac{1}{q'}}
=\| m\|_{q}(|\hat{f}|^{q'}*|\hat{g}|^{q'}*|\hat{h}|^{q'})(0)^{\frac{1}{q'}}
$$
$$
\leq \| m\|_{q}\| |\hat{f}|^{q'}*|\hat{g}|^{q'}*|\hat{h}|^{q'}\|_{\infty }^{\frac{1}{q'}}
\leq \| m\|_{q}
\Big( \| |\hat{f}|^{q'}\|_{\tilde{p}_{1}'}
\| |\hat{g}|^{q'}\|_{\tilde{p}_{2}'}
\| |\hat{h}|^{q'}\|_{\tilde{p}_{3}}\Big)^{\frac{1}{q'}}
$$
$$
=\| m\|_{q}
\| \hat{f}\|_{\tilde{p}_{1}'q'}\| \hat{g}\|_{\tilde{p}_{2}'q'}
\| \hat{h}\|_{\tilde{p}_{3}q'}
=\| m\|_{q}
\| \hat{f}\|_{p_{1}'}\| \hat{g}\|_{p_{2}'}\| \hat{h}\|_{p_{3}}
\leq \| m\|_{q}
\| f\|_{p_{1}}\| g\|_{p_{2}}\| h\|_{p_{3}'}
$$

\begin{nota}
Though $K\in L^{p}$ for some $1<p\leq 2$ none of the functions
$m(\xi ,\eta )=\hat{K}(\alpha \xi +\beta \eta )\in L^{q}(\bbbr^{2})$ for
$1\leq q\leq 4$. So, this result is neither a generalization nor a
particularization of proposition \ref{youngbilineal}.
\end{nota}

\bibliographystyle{amsplain}

\end{document}